\documentclass[11pt]{amsart}
\usepackage{amsfonts, amsmath, amssymb}
\newtheorem{theorem} {Theorem}
\newtheorem{lemma} [theorem]{Lemma} 

\newtheorem{cor} [theorem]{Corollary}
 
\newtheorem{prop}[theorem]{Proposition}
\newtheorem{rem}{Remark}

\newtheorem{fact}[theorem]{Fact}

\newcommand{\C}{{\mathbb C}}
\newcommand{\Q}{{\mathbb Q}}
\newcommand{\Z}{{\mathbb Z}}

\newcommand{\OO}{{\mathcal O}}
\newcommand{\CC}{{\mathbb C}}
\newcommand{\AAA}{{\mathcal A}}
\newcommand{\pp}{{\mathfrak p}}
\newcommand{\PP}{{\mathfrak P}}
\newcommand{\QQ}{{\mathfrak Q}}

\newcommand{\F}{{\mathbb F}}

\newcommand{\Gal}{\operatorname{Gal}}   

\newcommand{\mmod}{\operatorname{mod }}

\newcommand{\isom}{\cong}

\newcommand{\NN}{{\operatorname{\bf{N}}}}
\newcommand{\End}{\operatorname{End}}
\newcommand{\Ker}{\operatorname{Ker}}
\newcommand{\pibar}{\overline{\pi}}
\newcommand{\comment}[1]{}
\begin{document}
\title[A CRT algorithm for constructing genus 2 curves]
{A CRT algorithm for constructing genus 2 curves over finite fields}
\author{Kirsten Eisentr\"ager and Kristin Lauter} \address{Department
  of Mathematics, University of Michigan, Ann
  Arbor, MI 48109, USA.}  \email{eisentra@umich.edu} \address{Microsoft
  Research, One Microsoft Way, Redmond, WA 98052, USA.}
\email{klauter@microsoft.com}
\begin{abstract}
  We present a new method for constructing genus $2$ curves over a
  finite field $\F_n$ with a given number of points on its Jacobian.
  This method has important applications in cryptography, where groups
  of prime order are used as the basis for discrete-log based
  cryptosystems.  Our algorithm provides an alternative to the
  traditional CM method for constructing genus 2 curves. For a quartic
  CM field $K$ with primitive CM type, we compute the Igusa class
  polynomials modulo $p$ for certain small primes $p$ and then use the
  Chinese remainder theorem (CRT) and a bound on the denominators to
  construct the class polynomials. We also provide an algorithm for
  determining endomorphism rings of ordinary Jacobians of genus 2
  curves over finite fields, generalizing the work of Kohel for
  elliptic curves.
\end{abstract}
\thanks{The first author was partially supported by the National
  Science Foundation under agreement No.\ DMS-0111298 and by a
  National Science Foundation postdoctoral fellowship.}  
\maketitle

\section{Introduction}
In cryptography, some public key protocols for secure key exchange and
digital signatures are based on the difficulty of the discrete
logarithm problem in the underlying group.  In that setting, groups
such as the group of points on an elliptic curve or the group of
points on the Jacobian of a genus 2 hyperelliptic curve over a finite
field may be used.  The security of the system depends on the the
largest prime factor of the group order, and thus it is desirable to
be able to construct curves such that the resulting group order is
prime. This paper presents an alternative to the CM (Complex
Multiplication) algorithm for generating a genus 2 curve over a finite
field with a known number of points on its Jacobian.

The CM algorithm for genus $2$ is analogous to the Atkin-Morain CM
algorithm for elliptic curves proposed in the context of primality
testing (\cite{AM}).  Whereas the Atkin-Morain algorithm generates
the Hilbert class polynomial of an imaginary quadratic field $K$ by
evaluating the modular $j$-invariants of all elliptic curves with CM
by $K$, the genus $2$ algorithm generates what we will refer to as the
{\em Igusa class polynomials} of a quartic CM field $K$ by evaluating
the modular
invariants of all the abelian varieties of dimension $2$ with CM by
$K$.  Just as the $j$-invariant of an elliptic curve can be calculated
in two ways, either as the value of a modular function on a lattice
defining the elliptic curve as a complex torus over $\C$ or directly
from the coefficients of the equation defining the elliptic curve, the
triple of Igusa invariants (\cite{Igusa1,Igusa3}) of a genus $2$ curve
can also be calculated in two different ways.  Using classical
invariant theory over a field of characteristic zero, Clebsch defined
the triple of invariants of a binary sextic $f$ defining a genus $2$
curve $y^2=f(x)$.  Bolza showed how those invariants could also be
expressed in terms of theta functions on the period matrix associated
to the Jacobian variety and its canonical polarization over $\C$.
Igusa showed how these invariants could be extended to work in
arbitrary characteristic~(\cite[p.\ 848]{Igusa2}, see
also~\cite[Section 5.2]{GL}), and so the invariants are often referred
to as Igusa or Clebsch-Bolza-Igusa invariants.

To recover the equation of a genus $2$ curve given its invariants,
Mestre gave an algorithm which works in most cases, and involves
possibly passing to an extension of the field of definition of the
invariants (\cite{Mestre}).  The CM algorithm for genus $2$ curves
takes as input a quartic CM field $K$ and outputs the Igusa class
polynomials with coefficients in $\Q$ and if desired, a suitable prime
$p$ and a genus $2$ curve over $\F_p$ whose Jacobian has CM by $K$.
The CM algorithm has been implemented by Spallek~(\cite{Spallek}), van
Wamelen~(\cite{vanWamelen}), Weng~(\cite{Weng}),
Rodriguez-Villegas~(\cite{RV}), and Cohn-Lauter~(\cite{CL}).  This
method requires increasingly large amounts of precision of accuracy to
obtain the theta values necessary to form the class polynomials.  The
running time of the CM algorithm has not yet been analyzed due to the
fact that no bound on the denominators of the coefficients of the
Igusa class polynomials was known prior to the work of~\cite{GL}.

The idea of the algorithm we present here is to calculate the Igusa
class polynomials of a quartic CM field in a different way than the CM
algorithm does.  Our method generalizes the algorithm for finding the
Hilbert class polynomial given in \cite{ALV} to the genus 2 situation.
Given a quartic CM field $K$ with primitive CM type, for each small
prime $p$ in a certain set we determine the Igusa class polynomial
modulo $p$ by finding all triples of invariants modulo $p$ for which
the corresponding genus $2$ curve has CM by $K$.  The Igusa class
polynomial is then found using the Chinese Remainder Theorem (or the
explicit CRT as in~\cite{ALV}) and a bound on the denominators of the
coefficients.

Several difficulties arise in the genus $2$ situation which are absent
in the elliptic curve case. In this paper we resolve the following
issues: the field of definition of a CM abelian variety, necessary
conditions on the small primes for the algorithm to succeed, and the
computation of the endomorphism ring of the Jacobian of a genus $2$
curve in the ordinary case.  Our algorithm for computing endomorphism
rings of Jacobians of genus 2 curves over finite fields generalizes
the work of Kohel~\cite{Kohel} for elliptic curves.

\subsection{Statement of the Theorem} We will refer to a quartic CM field
$K$ with primitive CM type as a {\it primitive quartic CM field}.
Given a primitive quartic CM field $K$, let $\mathcal{A}$ be a system
of representatives for the set of isomorphism classes of principally
polarized abelian varieties over $\CC$ having complex multiplication
by $\OO_K$. For each abelian variety $A\in \mathcal{A}$ let
$(j_{1}(A),j_{2}(A),j_{3}(A))$ be the absolute Igusa invariants of
$A$.  Then the {\it Igusa class polynomials} $H_i$, for $i=1,2,3$, are
defined to be
\[
H_i := \prod_{A\in \mathcal{A}} (X- j_{i}(A)).
\]
It is known~(\cite{Shimura}) that roots of these polynomials generate
unramified abelian extensions of the reflex field of $K$.  It is also
known that Igusa class polynomials can be used to generate genus 2
curves with CM by $K$, and thus with a given zeta function over a
suitable prime field ({\it cf.\ }Section~\ref{S-curves}).  In this
paper we prove the following theorem.
\begin{theorem} \label{main} Given a quartic CM field $K$ with
  primitive CM type, the following algorithm finds the Igusa class
  polynomials of $K$:

{\bf (1)} Produce a collection $S$ of small rational primes $p \in S$
satisfying:

a. $p$ splits completely in $K$ and splits completely into principal ideals
in $K^*$, the reflex of $K$.  

b. Let $B$ be the set of all primes of bad reduction for the
genus 2 curves with CM by $K$. Then $S \cap B = \emptyset$.

c. $\prod_{p \in S} p >c$, where $c$ is a constant determined in
Theorem~\ref{CRTstep}.

{\bf (2)} Form the class polynomials $H_1$, $H_2$,
$H_3$ modulo $p$ for each $p \in S$.
Let $H_{i,p}(X):= H_i(X) \mod p$. Then
 \begin{equation*}
H_{i,p}(X)= \prod_{C \in T_p}(X-j_i(C)),\end{equation*}
where $T_p$ is the collection of $\overline{\F}_p$-isomorphism
classes of genus 2 curves over $\F_p$ whose Jacobian has
endomorphism ring isomorphic to $\OO_K$.

{\bf (3)} {\it Chinese Remainder Step.} Form $H_i(X)$
from $\{H_{i,p}\}_{p  \in S}$ $(i=1,2,3)$.
\end{theorem}

\begin{rem} Condition 1(a) is enough to insure that $p$ solves a
relative norm equation in $K/K_0$, $\pi \pibar = p$, $\pi$ a Weil
number ({\it cf.\ }Proposition~\ref{P-grouporder} below).  
\end{rem}
\begin{rem}
By~\cite{GL}, the primes in the set $B$ and in the denominators
of the class polynomials are bounded effectively by a quantity related
to the discriminant of $K$. Furthermore, it follows
from~\cite[Theorems 1 and 2]{Goren} and the discussion
in~\cite[Section 4.1]{GL} that condition 1(b) is implied by condition 1(a).
\end{rem}
\begin{rem}
  It follows from the Cebotarev density theorem that the density of
  the primes in the set $S$ is inversely proportional to the class
  number of $K$ in the case that $K$ is Galois cyclic. In the
  non-Galois case, the density is inversely proportional to the degree
  of the normal closure of the composite of $K$ with the Hilbert class
  field of the reflex of $K$.
\end{rem}

Our algorithm in the present form is not efficient, and we make no
claims about the running time.  A complete implementation of our
algorithm is now available in~\cite{FL}, along with new efficient
probabilistic algorithms for computing endomorphism rings.  Our
algorithm has the advantage that it does not require exponentially
large amounts of precision of computation.  It was recently brought to
our attention that the paper~\cite{Chao} proposes a similar algorithm,
but they give no proof of the validity of the approach.  Indeed, they
fail to impose the conditions necessary to make the algorithm correct
and include many unclear statements.  Also, while revising this paper,
a $p$-adic approach to generating genus 2 curves was given
in~\cite{GHKRW}.  No comparison has yet been made between the
different available approaches.

In Section~\ref{S-curves} we show how Theorem~\ref{main} can be used
to generate genus 2 curves with a given zeta function. The proof of
Theorem~\ref{main} is given in Section~\ref{correctness}.
Implementation details for the algorithm are given in
Section~\ref{S-implement}.  In Section~\ref{S-Kohel} we show how to
determine the endomorphism ring of an ordinary Jacobian of a genus 2
curve.  Section~\ref{S-example} gives an example of the computation of
a class polynomial modulo a small prime.

\

\noindent {\bf Acknowledgments.}  The authors thank E.\ Goren, E.\
Howe, K.\ Kedlaya, J-P.\ Serre, P.\ Stevenhagen, and T.\ Yang for
helpful discussions. The authors also thank D.\ Freeman and the
referee for valuable comments to improve the paper.

\section{Notation}\label{notation}  Throughout this paper, $C$ denotes
a smooth, projective, absolutely irreducible curve, and $J=J(C)$ will
be its Jacobian variety with identity element $\mathbf{O}$. The field
$K$ is always assumed to be a primitive quartic CM field, $K \neq \Q(\zeta_5)$, with ring
of integers $\OO_K$.  The real quadratic subfield of $K$ is denoted by
$K_0$, and a generator for the Galois group $\Gal(K/K_0)$ is denoted
by a bar, $\omega \mapsto \bar\omega$.  We will write $K^{*}$ for the reflex
of the quartic CM field $K$.  For $i=1,2,3$ we let $H_i(X)$ be the
Igusa class polynomials of $K$, and for a prime $p \in S$ we let
$H_{i,p}:=H_i\mod p$.  For a field $F$, $\overline{F}$ will denote an
algebraic closure of $F$.  We say that {\it $C$ has CM by $K$} if the
endomorphism ring of $J(C)$ is isomorphic to the full ring of integers
$\OO_K$.

\section{Generating genus 2 curves with 
a given zeta function} \label{S-curves}

Our algorithm solves the following problem under certain
conditions.

\noindent
{\bf Problem:} Given $(n, N_1, N_2)$, find a genus 2 curve $C$ over
the prime field $\F_n$ such that $\#C(\F_n)=N_1$ and
$\#C(\F_{n^2})=N_2$. 

Given $(n, \, N_1, \, N_2)$, it is straightforward to find $K$, the
quartic CM field such that the curve $C$ has CM by $K$, by finding the
quartic polynomial satisfied by Frobenius. Write $N_1=n+1-s_1,$ and
$N_2 = n^2+1+2s_2-s_1^2,$ and solve for $s_1$ and $s_2$. Then $K$ is
generated over $\Q$ by the polynomial $t^4-s_1t^3+s_2t^2-ns_1t+n^2$.

\noindent
{\bf Restrictions:}
If $s_2$ is prime to $n$, then the Jacobian is ordinary~(\cite[p.\ 
2366]{Howe}). Assume that $(s_2,n)=1$.  We also restrict to primitive
CM fields $K$. If $K$ is a quartic CM field, then $K$ is not primitive
iff $K/\Q$ is Galois and biquadratic ($\Gal(K/\Q)=V_4$) (\cite[p.\ 
64]{Shimura}).  In the example in Section~\ref{S-example}, $K$ is given
in the form $K=(i\sqrt{a+b\sqrt d})$, with $a,b,d \in \Z$ and $d$ and
$(a,b)$ square free. In this form the condition is easy to check: $K$
is primitive iff $a^2 -b^2d \neq k^2$ for some integer $k$ (\cite[p.\ 
135]{KaWa}).  Assume further that $K$ does not contain a cyclotomic
field.

\noindent
{\bf Solution:}
Given a triple $(n, \, N_1, \, N_2)$ satisfying the above
restrictions, one can generate a curve $C$ over $\F_n$ with the
associated zeta function as follows.  Compute $K$ and its Igusa class
polynomials $H_1, \, H_2, \, H_3$ using Theorem~\ref{main}.  From a
triple of roots modulo $n$ of $H_1, \, H_2, \, H_3$, construct a genus
2 curve over $\F_n$ using the combined algorithms of Mestre
(\cite{Mestre}) and Cardona-Quer (\cite{CQ}).  To match triples of
roots, in practice one can test whether the curve generated has the
correct zeta function by checking the number of points on the Jacobian
of the curve.  A curve $C$ with the correct zeta function will have
$\#J(C)(\F_n)=N=(N_1^2+N_2)/2 - n$.  If the curve does not have the
required number of points on the Jacobian, a twist of the curve may be
used.  In the case where 4 group orders are possible for the pair
$(n,K)$ ({\it cf.  }Section~\ref{S-grouporder}), a different triple of
invariants may be tried until the desired group order is obtained.

\section{Proof of Theorem~\ref{main}}\label{correctness}
Given a primitive quartic CM field $K$, let $\AAA$ be a system of
representatives of the isomorphism classes of simple principally polarized
abelian surfaces over $\C$ with CM by $K$. Each element of $\AAA$ has
a field of definition $k$ which is a finite extension of $\Q$
(\cite[Prop.\ 26, p.\ 96]{Shimura}).  For any prime $p \in S$
satisfying the conditions of Theorem~\ref{main}, the set $T_p$ was
defined in Step 2 of Theorem~\ref{main} as the collection of
$\overline{\F}_p$-isomorphism classes of genus 2 curves over $\F_p$
with an isomorphism of $\OO_K$ with $\End(J(C))$.  We claim that we have a
bijective correspondence between $\AAA$ and $T_p$. Moreover, we claim
that reducing the Igusa invariants gives the Igusa invariants of the
reduction.  Taken together, these can be stated in the form of the
following theorem:
\begin{theorem}\label{Step2proof}
  Let $K$ be a primitive quartic CM field and let $p \in S$ be a rational
  prime that satisfies the conditions of Theorem~\ref{main}. 
  Then $$H_{i,p}(X)= \prod_{C \in T_p}(X-j_i(C)),$$ where
  $H_{i,p}(X)$ and $T_p$ are defined as in Theorem~\ref{main}.   
\end{theorem}
\begin{proof}
  Let $A \in \AAA$ be a principally polarized abelian surface with CM
  by $K$, defined over a number field $k$.  Let $k_0$ be its field of
  moduli (see~\cite[p.\ 27]{Shimura} for the definition). By class
  field theory, $p$ splits completely into principal ideals in $K^*$
  if and only if $p$ splits completely in $H^*$, the maximal
  unramified abelian extension of $K^*$~(\cite[Corollary 5.25]{Cox}).
  The field of moduli $k_0$ is contained in $H^*$ (see~\cite[Main
  Theorem 1, p.\ 112]{Shimura}), but in general it is not true that
  $k=k_0$.  By a theorem of Shimura~(see \cite[Ex.\ 1, p.\
  525]{Shimura71}, see also~\cite[Proposition 2.1]{Goren}) if $K$ is a
  primitive quartic CM field, then $k$ is contained in $k_0$, so $A$
  is defined over $k_0$.
  
  Proposition 2.1 of \cite{Goren} also shows that $A$ has good
  reduction at any prime $\beta$ of $\OO_{H^*}$.  Let $A_p$ be the
  reduction of $A$ modulo a prime above $p$.  Then because $p$ splits
  completely in the Galois closure of $K$, $A_p$ is
  ordinary~(\cite[Theorems 1 and 2]{Goren}) and because $p$ splits
  completely into principal ideals in $K^*$, $A_p$ is defined over
  $\F_p$.  By condition 1(b) of Theorem~\ref{main}, $A_p$ is the
  Jacobian of a genus 2 curve $C$ over $\F_p$~(\cite{OortUeno}).  Then
  $C$ is an element of $T_p$.
  
  We must show that this correspondence is one-to-one and onto.  To
  show that it is one-to-one, we can generalize the argument in
  \cite[Theorem 13, p.\ 183]{Lang}. Let $A,B \in \AAA$, and for $p \in
  S$ let $A_p$ and $B_p$ be the reductions of $A$ and $B$ as above.
  Assume that $A_p$ and $B_p$ are isomorphic over $\overline{\F}_p$,
  and let $\varepsilon: B_p \to A_p$ be an isomorphism. The varieties
  $A$ and $B$ both have CM by $K$, hence there exists an isogeny $\lambda:A
  \to B$ (\cite[Corollary, p.\ 41]{Shimura}) giving rise to a reduced
  isogeny $\lambda_p: A_p \to B_p$.  Since the endomorphism ring of $A$ is
  preserved under the reduction map, there exists $\alpha \in \End(A)$ such
  that the reduction $\alpha_p$ satisfies $\alpha_p = \varepsilon \circ \lambda_p$.
  Let $C$ be the image of the map $\lambda \times \alpha: A \times A \to B \times A $. With a
  similar argument as in \cite[p.\ 184]{Lang}, one can then show that
  $C$ is the graph of an isomorphism between $A$ and $B$. Similarly,
  if there is an isomorphism of the principal polarizations on $A_p$
  and $B_p$ then this isomorphism lifts to an isomorphism of the
  polarizations on $A$ and $B$.  This shows that the correspondence is
  one-to-one.

  The correspondence is onto because, given a genus 2 curve $C$ over
  $\F_p$ with CM by $K$ representing a class of $T_p$, its Jacobian
  $J(C)$ is ordinary and so it can be lifted, along with its
  endomorphism ring and its polarization, to its 
``Serre-Tate canonical lift'', $A$, defined
  over the Witt vectors $W(\F_p)=\Z_p$ (\cite[Theorem 3.3,
  p.\ 172]{Messing}). Let $L$ be the field generated over $\Q$ by all the
  coefficients of the equations defining $A$. Then $A$ is defined over
  $L$ and since $L$ has finite transcendence degree over $\Q$, we can
  embed it into $\C$.  So we can lift $J(C)$ to an abelian variety
  with CM by $K$ defined over $\C$.
  
  By assumption 1(b) of Theorem~\ref{main}, no prime above $p \in S$
  is a prime of bad reduction for a genus 2 curve with CM by $K$, so
  by~\cite[Cor 5.1.2]{GL}, $p \in S$ is coprime to the denominators of
  the class polynomials $H_i(X)$.  We claim that reducing the
  coefficients of $H_i$ modulo $p$ gives the same result as taking the
  polynomial whose roots are the absolute Igusa invariants of the
  curves over $\F_p$ with Jacobians equal to the reductions modulo a
  prime above $p$ of the abelian varieties $A$ representing the
  classes of $\AAA$.  Since the absolute Igusa invariants are rational
  functions in the coefficients of the curve, the order of computation
  of the invariants and reduction modulo a prime can be reversed as
  long as the primes in the denominator are avoided and an appropriate
  model for the curve is chosen.
\end{proof}

\begin{theorem}\label{CRTstep}
  Suppose the factorization of the denominators of the
  Igusa class polynomials is known.  Let $\nu$ be the
  largest absolute value of the coefficients of the $H_i$, and let $\lambda$ be the least common multiple of the
  denominators of the coefficients of the $H_i$ $(i=1,2,3)$.  Let $S$
  be a set of rational primes such that $S \cap B = \emptyset$ and
  $\prod_{p\in S} p >c$, where $c=2\lambda \cdot \nu$. Then the Chinese
  Remainder Theorem can be used to compute the class polynomials
  $H_i(X) \in \Q[X]$ from the collection $\{H_{i,p}\}_{p \in S}$,
  $i=1,2,3$.
\end{theorem}
\begin{proof}
  By assumption $\lambda$ is prime to all $p \in S$. The polynomials
  $$F_i(X):=\lambda \cdot H_i(X) \,\, i=1,2,3$$
  have integer coefficients.  For each $p\in S$ let
  $$F_{i,p}:= F_i \, (\mmod p) = \lambda \cdot H_{i,p} \, (\mmod p).$$
  Apply the Chinese Remainder Theorem to the collection
  $\{F_{i,p}\}_{p\in S}$ to obtain a polynomial which is congruent to
  $F_i \in \Z[X]$ modulo the product $\prod_{p \in S} p$. Since $c$
  was taken to be twice $\lambda$ times the largest absolute value of the
  coefficients, we have found $F_i$, and so $H_i = \lambda^{-1}\cdot
  F_i$.
\end{proof}
\begin{rem}
  It was proved in~\cite{GL} that the primes dividing the denominators
  are bounded effectively in terms of the field $K$ by a quantity
  related to the discriminant.  The power to which each prime in the
  denominator appears has also been bounded in recent work of Goren,
  and so we can conclude that we have a bound on the denominators of
  the class polynomials.
\end{rem}
\begin{proof}[Proof of Theorem~\ref{main}]
  The proof of Theorem~\ref{main} now follows immediately from
  Theorem~\ref{Step2proof} and Theorem~\ref{CRTstep}.
\end{proof}
\section{Implementation} \label{S-implement}
\subsection{The possible group orders for each p}
\label{S-grouporder}
Suppose that $C$ is a genus 2 curve defined over $\F_p$ with CM by
$K$. To find all possible group orders for $J(C)(\F_p)$, let $\pi \in
O_K$ correspond to the Frobenius endomorphism of $C$.  Since the
Frobenius satisfies $\pi\pibar = p$, it follows that the relative norm
of $\pi$ is p, i.e.\ $\NN_{K/K_0}(\pi)=p$, and hence
$\NN(\pi)=\NN_{K/\Q}(\pi)=p^2$.  So if $K$ is fixed, primes $p$ for
which there exist genus 2 curves modulo $p$ with CM by $K$ are primes
for which there are solutions to the relative norm equation:
$\NN_{K/K_0}(\pi)=p$.  The following proposition gives the number of possible group
orders in each case.  It overlaps with~\cite[Thm 4.1]{Weng2}, but our statement, assumptions,
and proof are all slightly different, and we use the details of this proof in our algorithm, so we
include it here. Note that, as pointed out in~\cite{Weng2}, it is not known whether two of the four possible group orders could coincide in the non-Galois case.

\begin{prop} \label{P-grouporder}
Fix a primitive quartic CM field $K$, and a rational prime $p$ unramified in $K$.
Assume that $K \neq \Q(\zeta_5)$, so that the only roots of unity in $K$ are $\{\pm 1 \}$. Then

(A) There are either 0, 2 or 4 possibilities for the group order
  $\#J(C)(\F_p)$ of curves $C$ with CM by $K$.

(B) Under the additional assumption that $p$ splits completely into 
principal ideals in $K^*$ and splits completely in $K$, there are always 
2 possible group orders in the cyclic case and 4 possible group orders 
in the non-Galois case.
\end{prop}
\begin{proof}
  We consider all possible decompositions of the prime $p$ in $K$. \\[0,1cm]
  {\bf Case 1:} There exists a prime ideal $\pp$ of $K_0$ above $p$
  that does not split in $K$. In this case there is no
  solution to the relative norm equation.\\[0,1cm]
  {\bf Case 2:} The rational prime $p$ is inert in $K_0/\Q$, and the
  prime $\pp$ of $K_0$ above $p$ splits in $K$ with $\PP_1|\pp$ and
  $\PP_2|\pp$.  We have $\overline{\PP_1}=\PP_2$. In this case there
  are two ideals of norm $p^2$, $\PP_1$ and $\PP_2$. If $\PP_1$ is not
  principal, then there are no solutions to the norm equation. If
  $\PP_1$ is principal with generator $\pi$, then $\PP_2 = (\pibar)$,
  and $\pi \pibar = p$.  The elements $\pi$ and $\pibar$ are Galois
  conjugates, so by Honda-Tate $\pi$ and $-\pi$ give rise to all
  possible group orders.  Let $\pi_1:=\pi$, and let $\pi_2, \dots,
  \pi_4$ be its conjugates over $\Q$.  Then $m_1 =
  \prod_{i=1}^{4}(1-\pi_i)$ and $m_2 = \prod_{i=1}^{4}(1-(-\pi_i))$
  are the 2 possible group orders for the Jacobian. \\[0,1cm]
  {\bf Case 3:} $p$ splits completely in $K/\Q$, with
  $\PP_1,\dots,\PP_4$ lying above $p$ and with
  $\overline{\PP_1}=\PP_2$, and $\overline{\PP_3}=\PP_4$.  Then $\PP
  := \PP_1 \PP_3$ , $\QQ := \PP_1 \PP_4$, and $\overline{\PP}$ and
  $\overline{\QQ}$ are the only ideals with relative norm $p$.

{\bf Subcase (a)} If $K/\Q$ is Galois, then the Galois group is
cyclic, since we assumed that $K$ was a primitive CM field (\cite[p.\ 
65]{Shimura}).  Let $\sigma$ be a generator of $\Gal(K/\Q)$. Then
w.l.o.g.  $\PP_2=\PP_1^{\sigma^2}, \PP_3=\PP_1^{\sigma}$, and
$\PP_4=\PP_1^{\sigma^3}$.  Thus $\PP = \PP_1 \PP_1^{\sigma}=
(\PP_1\PP_1^{\sigma^3})^{\sigma}=\QQ^{\sigma}$, so if $\PP$ is
principal, so is $\QQ$, and their generators, $\omega$ and
$\omega^{\sigma}$ give rise to isogenous curves. Hence if $\PP$ is
principal, then there are two possible group orders as before, and if
it is not principal, then the relative norm equation has no solution.

{\bf Subcase (b)} If $K/\Q$ is not Galois, then the Galois group of
its splitting field is the dihedral group $D_4$ (\cite[p.\ 
65]{Shimura}). In this case $\PP$ and $\QQ$ are not Galois conjugates.
So if both $\PP$ and $\QQ$ are principal, then there are 4 possible
group orders, if only one of them is principal, then there are 2
possible group orders, and otherwise there are no solutions to the
relative norm equation.

Statement (A) follows from the 3 cases considered above.  Statement
(B) concerns Case 3.  If $K$ is Galois, then $K=K^{*}$ and the
additional assumptions imply that $\PP$ is principal, and then there
are 2 possible group orders.  If $K$ is not Galois, let $L$ be the
Galois closure with dihedral Galois group $\Gal(L/\Q) = \langle
\tau,\sigma:$ $\tau^2, \sigma^4, \tau \sigma \tau \sigma \rangle$ such that
$K$ is the fixed field of $\tau$ and the CM type is $\{1,\sigma \}$.
Then $\sigma^2$ is complex conjugation.  According to~\cite[Theorem
2]{Goren}, a rational prime $p$ that splits completely in $L$ with
$\mathcal{P}:= p\mathcal{O}_L$ decomposes as follows in $K$ and
$K^*$:
$$p\mathcal{O}_K = \PP_1\PP_2\PP_3\PP_4=(\mathcal{P}\mathcal{P}^{\tau})
(\mathcal{P}^{\sigma^2}\mathcal{P}^{\tau \sigma^2})
(\mathcal{P}^{\sigma}\mathcal{P}^{\tau \sigma})
(\mathcal{P}^{\sigma^3}\mathcal{P}^{\tau \sigma^3}),$$
$$p\mathcal{O}_{K^*} = \PP_1^*\PP_2^*\PP_3^*\PP_4^*= 
(\mathcal{P}\mathcal{P}^{\tau \sigma^3})
(\mathcal{P}^{\sigma^2}\mathcal{P}^{\tau \sigma})
(\mathcal{P}^{\sigma}\mathcal{P}^{\tau})
(\mathcal{P}^{\sigma^3}\mathcal{P}^{\tau \sigma^2}).$$
By assumption, $\PP_1^*$, $\PP_2^*$, $\PP_3^*$, $\PP_4^*$ are principal.
Thus both $\PP$ and $\QQ$ are principal
since
$\PP = \PP_1 \PP_3 = \PP_3^* (\PP_4^*)^\sigma,$
and
$\QQ = \PP_1 \PP_4 = \PP_1^* (\PP_1^*)^{\tau}.$
Thus there are 4 possible group orders when $K$ is not Galois.
\end{proof}

\subsection{Generating the collection of primes $S$} \label{S-primes}
In practice to generate a collection of primes belonging to $S$
there are several alternatives. One approach is to run
through small primes checking the splitting behavior in $K$ and $K^*$
using a computational number theory software package like PARI.  
A second approach is to generate solutions to the relative norm
equation directly as in~\cite[Section 8]{Weng}, then check each solution for the splitting in $K$
and $K^*$ and check for the other solution to the relative norm
equation in the case that $K$ is not Galois.  One advantage to this
approach is that it gives direct control over the index of
$\Z[\pi,\pibar]$ in $\OO_K$ in terms of the coefficients $c_i$ of
$\pi$, the solution to the relative norm equation ({\it cf.\ 
}Proposition~\ref{P-bound}).

\subsection{Computing Igusa class polynomials modulo p}\label{S-classpoly}
Let $p \in S$. To compute the Igusa class polynomials mod $p$ we must
find all $\overline{\F}_p$-isomorphism classes of genus 2 curves over
$\F_p$ whose Jacobian has CM by $K$.  This can be done as follows:

{\bf (1)} For each triple of Igusa invariants modulo $p$, generate a
genus 2 curve with those Igusa invariants using an implementation of
the Mestre-Cardona-Quer algorithm (\cite{Mestre}, \cite{CQ}).
  
  {\bf (2)} Let $N_p:=\{(n_1,m_1),(n_2,m_2),\dots,(n_r,m_r)\}$ be the
  set of possible group orders $(\#C(\F_p),\#J(C)(\F_p))$ of curves
  $C$ which have CM by $K$ as computed above in Section~\ref{S-grouporder}. 
  
  {\bf (3)} Collect all curves $C$ such that $(\#C(\F_p),\#J(C)(\F_p))
  \in N_p$ as follows: for each triple of invariants and a
  corresponding curve $C$, take a random point $Q$ on $J(C)$. Multiply
  $Q$ by $m_1,\dots,m_r$ and check if the identity element is obtained
  for some $r$. If not, then $C$ does not belong to $T_p$.  If a curve
  passes this test, then count the number of points on the curve and
  its Jacobian over $\F_p$ to check whether the Jacobian has the right
  isogeny type. This procedure obtains all curves in the desired
  isogeny class.  For each curve in the desired isogeny class, the
  endomorphism ring of the Jacobian contains the ring $\Z[\pi,\pibar]$
  and is contained in the ring $\OO_K$. The curve is included in the
  set $T_p$ only if $\End_{\F_p}(J(C))=\OO_K$. In the next section, we will
show how to test this property by computing the endomorphism ring $\End_{\F_p}(J(C))$.
\section{Computing endomorphism rings of genus 2 curves}
\label{S-Kohel}
\subsection{The index of $\Z[\pi,\pibar]$ in $\OO_K$} \label{index}

For a prime $p$ and a Frobenius element $\pi \in \OO_K$,
the smaller the index of $\Z[\pi,\pibar]$ in $\OO_K$,
the less work it takes to compute the endomorphism ring.  For
example, if the index is 1, then we can determine whether $C \in T_p$
just from counting points on $C$ and its Jacobian.  
Proposition~\ref{P-bound} gives a bound for the index of 
$\Z[\pi,\pibar]$ in $\OO_K$.
\begin{prop} \label{P-bound}
  Let $K:= \Q(\eta)$ be a quartic CM field, where $\eta =
  i\sqrt{a+b\sqrt{d}}$ with $a,b,d \in \Z$ and $d$ and $(a,b)$ square
  free. Let $\OO_K$ be its ring of integers.  Assume for simplicity
  that the Frobenius endomorphism of $C$ is of the form $\pi := c_1 +
  c_2 \sqrt{d} + (c_3 + c_4\sqrt{d})\eta$ with $c_1,\dots, c_4 \in
  \Z$, that $a^2 -b^2 d$ is square free and that the real quadratic subfield 
$K_0$ has class number 1. If $d  \equiv 2,3 \mmod 4$, then $[\OO_K:\Z[\pi,\pibar]]$ divides $8c_2(c_3^2  -c_4^2d)$. If $d \equiv 1 \mmod 4$, then $[\OO_K:\Z[\pi,\pibar]]$ divides $16c_2(c_3^2 -c_4^2d)$.
\end{prop}
\begin{proof}
We have
\begin{equation}
\pi + \pibar - 2 c_1= 2 c_2 \sqrt{d},
\end{equation}
\begin{gather}
[2c_2c_3-c_4(\pi + \pibar-2c_1 )](\pi - \pibar)= 4 c_2 (c_3^2 - c_4^2
d)\eta,\\
(c_3-c_4\sqrt{d})(\pi - \pibar)= 2(c_3^2 -c_4^2d)\eta.
\end{gather}
So $\Z[2c_2\sqrt{d},4 c_2 (c_3^2 - c_4^2 d)\eta ] \subseteq
\Z[\pi,\pibar]$.  Since $K_0$ has class number 1, we have a relative integral
basis of $\OO_K$ over $\OO_{K_0}$. We can choose a relative basis of
the form $\{1, \kappa\}$, and by \cite{SpeWi}, in the case that $d
\equiv 2,3 \mmod 4$, $\kappa$ is either
\[1. \;\eta/2\;\;\;\;\;\;
2.\; (1 + \eta)/2\;\;\;\;\;\; 3.\; (\sqrt{d}+\eta)/2\;\;\;\;\;\;\;
4.\; (1+\sqrt{d}+ \eta)/2.\] In each case the index of
$\Z[\sqrt{d},\eta]$ in $\OO_K$ is 2.
For $d \equiv 1 \mmod 4$, $\kappa$ is either
\[
5.\; (1 + \sqrt{d}+2 \eta)/4 \;\;\;\;\; 
6.\; (-1 + \sqrt{d} + \eta)/4 \;\;\;\;\;
7.\; (-b+ \sqrt{d} + 2 \eta)/4.
\]
Here, in each case the index of $\Z[\sqrt{d},\eta]$ in $\OO_K$ is 4.
We have
\begin{equation*}
\Z[\pi,\pibar]\subseteq \Z[\pi,\pibar,\sqrt{d}]\subseteq\Z[\sqrt{d},
\eta]\subseteq\OO_K,
\end{equation*}
with
$[\Z[\pi,\pibar,\sqrt{d}]:\Z[\pi,\pibar]]$ dividing $2c_2$ and
$ [\Z[\sqrt{d},\eta]:\Z[\pi,\pibar,\sqrt{d}]]$ dividing $2(c_3^2
-c_4^2d).$
If $d \equiv 2, \, 3 \, \mmod 4$, then $[\OO_K:\Z[\sqrt{d},\eta]]=2$, and hence the index 
$[\OO_K:\Z[\pi,\pibar]]$ divides $8c_2(c_3^2 -c_4^2d)$.
If $d \equiv 1 \mmod 4$, then  $[\OO_K:\Z[\sqrt{d},\eta]]=4$, and hence
$[\OO_K:\Z[\pi,\pibar]]$ divides $16c_2(c_3^2 -c_4^2d)$.  Since the index is a positive integer, it is thus also
bounded by these quantities.
\end{proof}
So if we want to minimize the index $[\OO_K:\Z[\pi,\pibar]]$ then we
have to minimize $c_2(c_3^2 -c_4^2d)$.  When $a^2 -b^2 d$ is not
square free the representation of the ring of integers can become more
complicated (\cite{SpeWi}), but the term we need to minimize is still
$c_2(c_3^2 -c_4^2d)$.  Using the relative basis of $\OO_K$ over
$\OO_{K_0}$ we can also determine which denominators can occur in the
coefficients $c_i$ of the Frobenius endomorphism and generalize our
argument to the general case.

\subsection{Determining the index of $\End(J)$ in $\OO_K$}
We can summarize the necessary conditions to ensure that
$[\OO_K:\End(J)]=1$ as follows:
\begin{lemma}
  Under the conditions of Section~\ref{index}, to show that the
  endomorphism ring of a curve is the full ring of integers $\OO_K$,
  it is sufficient to test whether:
\begin{enumerate}
\item $\sqrt{d}$ is an endomorphism, where $2c_2\sqrt{d}= \pi + \pibar
  - 2c_1.$
\item $\eta$ is an endomorphism, where
  $$(4c_2(c_3^2-c_4^2d))\eta =
  (2c_2c_3-c_4(\pi+\pibar-2c_1))(\pi-\pibar).$$
  Here the $c_i$'s are
  the coefficients of $\pi$ written in the relative basis.
\item $\kappa$ is an endomorphism, where $\kappa$ is one of the 7
  possible elements listed in Section~\ref{index} in the case that
  $a^2-b^2d$ is square free.
\end{enumerate}
\end{lemma}
If any one of these conditions fails, we conclude that the
endomorphism ring of the curve is not the full ring of integers
$\OO_K$.  When $a^2 -b^2d$ is not square free then the relative
integral basis is listed in the table in~\cite[p.\ 186]{SpeWi}.  This
algorithm can also be modified to test whether the endomorphism ring
of the curve is some other subring of $\OO_K$ or to compute the
endomorphism ring exactly.

To test whether $\sqrt{d}$, $\eta$, and $\kappa$ are endomorphisms, we
express them as above as polynomials in $\pi$ and $\pibar$ with
integral denominators determined by the $c_i$. It will be proved in
Section~\ref{action} below that in each case it suffices to check
whether the numerator acts as zero on the $s$-torsion, where $s$ is
the denominator.
\subsection{Action on $s$-torsion}\label{action}
\begin{prop}\label{factor}
  Assume that $k$ is an algebraically closed field and that $A,B,C$
  are abelian varieties over $k$. Let $\beta: A\to B$, $\gamma
  : A \to C$ be two isogenies with $\beta$ separable and
  $\Ker(\beta) \subseteq \Ker(\gamma)$.  Then there is a homomorphism
  $\delta:B \to C$ such that $\delta \cdot \beta = \gamma$.
\end{prop}
\begin{proof}
  This proof follows the argument of Remark 7.12 in \cite[p.\ 
  37]{Milne}.  Since $\beta$ is separable, we can form the quotient
  abelian variety $A/\Ker(\beta).$ From the universal property of
  $A/\Ker(\beta)$ we have a regular map $A/\Ker \beta \to B$,
  which is again separable and bijective. Since $B$ is nonsingular,
  this implies that it is an isomorphism. Thus $B\isom A/\Ker
  (\beta)$. After identifying $B$ with $A/\Ker(\beta)$ and using the
  universal properties of quotients again we find that there is a
  unique regular map $\delta$ such that $\delta \cdot \beta = \gamma$.
  Moreover, $\delta$ is automatically a homomorphism because it maps
  $\mathbf{O}$ to $\mathbf{O}$.
\end{proof}
\begin{prop}\label{general}
  Let $k$ be an algebraically closed field and let $A$ be an abelian
  variety over $k$. Let $R:=\End_k A$. Let $s \in R$ be separable and
  let $A[s]=\{P\in A(k):sP=\mathbf{O}\}=\Ker(s)$. Then $A[s]$ is a
  faithful $R/Rs$-module.
\end{prop}
\begin{proof}
  Clearly, $A[s]$ is an $R/Rs$-module. We have to show that $A[s]$ is
  a faithful $R/Rs$-module; that is, any $r\in R$ with $r \cdot A[s]
  =0$ belongs to $Rs$.  Suppose $r$ is such that $r \cdot A[s]=0$.
  Since $s$ is separable, this implies that $r=ts$ for some
  endomorphism $t$ of $A$ by Proposition~\ref{factor} above applied
  with $A=B=C$, $\beta = s$ and $\gamma =r$. This implies that $r \in
  Rs$, which proves the claim.
\end{proof}
We will frequently use the following
\begin{cor}\label{divisible}
  Let $A,k$ be as in Proposition~\ref{general}. Let $n$ be a positive
  integer coprime to the characteristic of $k$.  Suppose that $\alpha
  : A\to A$ is an endomorphism, with $A[n] \subseteq
  Ker(\alpha)$, i.e. $\alpha$ acts as zero on the $n$-torsion. Then
  $\alpha = \beta \cdot n = n \cdot \beta$, for some endomorphism
  $\beta$, i.e. $\alpha$ is divisible by $n$ in $R=\End_k(A)$.
\end{cor}

\subsection{Computing the index using division polynomials}

In~\cite{Cantor}, Cantor finds recursive formulae for division
polynomials for hyperelliptic curves with one point at infinity,
$P_{\infty}$.  The $r$th division polynomials he defines are
$(\delta_r(X),\epsilon_r(X))$ such that
$(\delta_r(\frac{x-X}{4y^2}),\epsilon_r(\frac{x-X}{4y^2}))$ represents
$r \cdot (x,y)$, where $(x,y)$ is a point on the curve thought of as
the point $(x,y) - P_{\infty}$ on the Jacobian.  For a general point
on the Jacobian represented as $D = P_1 + P_2 - P_{\infty}$, we see
that $r D = 0$ iff $r P_1 = - r P_2$. If $P_1 = (x_1,y_1)$ and $P_2 =
(x_2,y_2)$, then we can write down a system of equations and an ideal,
$I_r$, defining the solutions to the system, where $I_r$ is an ideal
in $\F_p[x_1,x_2,y_1,y_2]$. Various ways of finding the ideal $I_r$
have been investigated, from Gr\"obner bases to resultant computations
(see~\cite{GH} and \cite{GS}).

The ideal $I_r$ can be used to test the action of endomorphisms on the
$r$-torsion.  For example, to check that $\pi^k$ (or any other
polynomial in $\pi$) acts like $a$ on the $r$-torsion, it suffices to
check that in $\F_p[x_1,x_2,y_1,y_2]$,
$$\pi^k(D) \equiv aD \mod I_r.$$

Even if the best method for computing the $I_r$ is not yet completely
well understood in practice, in theory this is likely the most
efficient way to compute the action of endomorphisms on $r$-torsion.

\subsection{Computing the index through direct computation of the action
of Frobenius on the torsion subgroups}

In practice, we used a computational number theory software package
like MAGMA to compute the group structure of $J(C)(\F_{p^k})$ for
small values of $k$. Using the generators of $J(C)(\F_{p^k})$ we then
explicitly computed the action of Frobenius on various torsion
subgroups to determine whether or not certain elements of the ring of
integers are endomorphisms.  An example will be given in the next
section. In the example we will use the following fact repeatedly:
\begin{fact}\label{torsion}
  Let $\gamma_k$ be a positive integer coprime to $p$. All
  the $\gamma_k$-torsion is defined over $\F_{p^k}$ if and only if
  $\frac{\pi^k-1}{\gamma_k}$ is an endomorphism.
\end{fact}
Fact~\ref{torsion} follows immediately from Corollary~\ref{divisible}.
Note that it is {\it not} true in general that the field of definition
of the $r$-torsion for all $m$ is enough to determine the endomorphism
ring.  We found examples of curves where the field of definition
of the $r$-torsion was the same for all $r$, but the endomorphism
rings were different because the action of Frobenius on the torsion
subgroups was different.  However, there are special cases where
checking the field of definition of the torsion is enough:
\begin{rem}
In the case where $\OO_K$ is generated by elements of the form
$\frac{\pi^k-1}{\gamma_k}$, for some collection of pairs of integers
$(k,\gamma_k)$, then equality of the endomorphism ring with $\OO_K$
can be checked simply by checking the field of definition of the
$\gamma_k$-torsion.
\end{rem}
\section{Example} \label{S-example}
Let $K:=\Q(i\sqrt{13-3\sqrt{13}})$. In this example we will find the
Igusa class polynomials of $K$ modulo $43$ by finding all genus 2
curves $C$ defined over $\F_{43}$ (up to isomorphism over the
algebraic closure of $\F_{43}$) such that $\End(J(C)) \isom \OO_K$,
where $\OO_K$ is the ring of integers of $K$.  Let $K^*$ be the reflex
of $K$. Since $a^2-b^2d=2^2\cdot 13$, the extension $K/\Q$ is cyclic
(\cite[p.\ 88]{KaWa}), and hence $K^{*}=K$ (\cite[p.\ 65]{Shimura}).
The real quadratic subfield of $K$ is $K_0 := \Q(\sqrt{13})$.  The
prime 43 splits completely in $K=K^{*}$.  The
class number of $K$ is 2, and so since $K$ is Galois, we expect two
classes of curves over $\F_{43}$ with CM by $K$.  Let $\eta := i
\sqrt{13 - 3\sqrt{13}}$.  The ring of integers of $K$ is
\[\OO_K = \Z + \frac{\sqrt{13}+1}{2} \Z + (\Z +
\frac{\sqrt{13}+1}{2}\Z)\,\eta.\] 
Let $\delta:=({1+\sqrt{13}})/{2}$. The prime $43$ factors in $K/K_0$ as:
\begin{gather*}
43 = \pi \cdot {\overline{\pi}} = (-3+2\cdot \delta +
(-2-\delta)\,\eta) \cdot (-3+2\cdot \delta) +
(2+\delta)\,\eta).
\end{gather*}
The characteristic polynomial of the Frobenius element corresponding
to $\pi$ is
\[
\psi(t)= 1849t^4 + 344t^3 + 50t^2 + 8t + 1.
\]

Let $C$ be a curve over $\F_{43}$ whose Frobenius is $\pm \pi$.  Then
the possibilities for $(\#C(\F_{43}),\#J(C)(\F_{43}))$ are $(52,2252)$
and $(36,1548)$.  Using MAGMA we found (up to isomorphism over
${\overline{\F}_{43}}$) 67 curves whose Frobenius is $\pm \pi$. 
However, not all 67 curves have endomorphism ring equal to the
full ring of integers.  To eliminate those with smaller endomorphism
ring, we first observe that
$$\frac{\pi^4-1}{12}=
-2+24\sqrt{13}+\frac{17}{2}\sqrt{13}i\sqrt{13-3\sqrt{13}}
+\frac{113}{2}i\sqrt{13-3\sqrt{13}} \in \OO_K.$$
Then Fact~\ref{torsion} implies that any curve whose endomorphism ring is
the full ring of integers must have the full $12$-torsion defined over
$\F_{43}^4$.  We can check that this eliminates all but 6 of the 67
curves.  The Igusa invariants of the remaining 6
curves are:
\[(3,24,36),
(4,29,28), (29,24,13), (20,21,29), (20,23,19), (36,21,6).\] We expect
only 2 curves over $\F_{43}$ (up to isomorphism) with CM by $K$.  To
eliminate the other 4 curves from this list, it is enough in this case
to check the action of Frobenius on the 4-torsion.  By
Corollary~\ref{divisible}, $\delta = \frac{\pi + \pibar +6}{4}$ is an
endomorphism of $J(C)$ if and only if $\pi + \pibar +6$ acts as zero
on the 4-torsion, or equivalently, $\pi + \pibar$ acts as
multiplication-by-2 on the 4-torsion.

Consider a curve $C$ with Igusa invariants $(20,23,19)$ given by the
equation $C: y^2 = 5x^6+ 21 x^5 + 36 x^4 + 7 x^3 + 29 x^2 + 32 x +10$
over $\F_{43}$.  All the 4-torsion is defined over a degree 4
extension, and we can use MAGMA to compute a basis for the 4-torsion
by computing the abelian group structure over the degree 4 extension.

We can then compute that the action of Frobenius on the 4-torsion is
given in terms of some basis by the matrix $F$, and the action of $\pibar$ is given given by
$V$:
\[
F=\left(
\begin{matrix}
 1 & 0 & 1 & 3\\
2 & 1 & 1&  0 \\
0 & 2 &3 &2\\
2 &2 &2 &3
\end{matrix}
\right)\;\;\;\;
V=\left(
\begin{matrix}
1 &0 &3& 1\\
2& 1& 3& 0\\
0& 2& 3& 2\\
2& 2& 2& 3
\end{matrix}
\right).
\]
From this it is easy to see that $\pi + \pibar = [2]$ on the
4-torsion, so $\delta$ is an endomorphism of $C$.  Performing the
identical computation on a curve $C$ with Igusa invariants
$(36,21,6)$, we find that $\delta$ is also an endomorphism for this
curve.  Doing the same calculation for the remaining 4 triples of
Igusa invariants $(3,24,36),(4,29,28),(29,24,13),(20,21,29)$, we see
that $\pi + \pibar =[2]$ does not hold on the 4-torsion in those
cases, so $\delta \notin \End(J(C))$ for any of the corresponding
curves.

It is easy to see in this case that $\delta \in \End(J(C))$ and
$\frac{\pi^4-1}{12} \in \End(J(C))$ is enough to conclude that
$\End(J(C))=\OO_K$.  Hence the two triples of invariants corresponding
to curves with CM by K are $(36,21,6)$ and $(20,23,19)$.
In conclusion, we have obtained the three Igusa class polynomials modulo $43$
with our method:
$$H_{1,43}(X) = X^2 + 30X + 32,$$
$$H_{2,43}(X) = X^2 + 42X + 10,$$
$$H_{3,43}(X) = X^2 + 18X + 28.$$
These indeed agree modulo $43$ with the class polynomials with rational 
coefficients computed by evaluating the quotients of Siegel modular forms
with 200 digits of precision
as computed by van Wamelen~(\cite{vanWamelen}):\\[0,2cm]
\begin{small}

  $H_1(X)= X^2 - \frac{9625430292534239443768093859336546624656066801331680515511924}{1224160503138337270992732796402545210705949947}X + \\
   \frac{17211893103548805144815938862454140808252633213039291208686119112918076788941674683411636004}{58670687646017062528338814934164161420328368922180746779053222569},$

  $H_2(X)= X^2 - \frac{3237631624959669936998571242515324335027260}{7973132502458523379282597629}X + \\
\frac{101869481833026643236326057638275086345512388711354393815337676100}{387742378329008606934824201506984053723129},$

  $H_3(X) = X^2 - \frac{2511631949170772694805531862232571975071932}{23919397507375570137847792887}X + \\
   \frac{83671593583457548222292142563905819629154823011540406083420061764}{3489681404961077462413417813562856483508161}.$
\end{small}


\begin{thebibliography}{CMKT00}

\bibitem[ALV04]{ALV}
Amod Agashe, Kristin Lauter, Ramarathnam Venkatesan.
\newblock {Constructing elliptic curves with a known number of points over a
  prime field}.
\newblock In {\em High Primes and Misdemeanours: lectures in honour of
  the 60th birthday of Hugh Cowie Williams}, Fields Institute Communications
  Series, Vol. 42, pp. 1--17, 2004.

\bibitem[AM93]{AM}
A.~O.~L. Atkin and F.~Morain.
\newblock Elliptic curves and primality proving.
\newblock {\em Math. Comp.}, 61(203):29--68, 1993.

\bibitem[Can94]{Cantor}
David~G. Cantor.
\newblock On the analogue of the division polynomials for hyperelliptic curves.
\newblock {\em J. Reine Angew. Math.}, 447:91--145, 1994.


\bibitem[CQ05]{CQ} Gabriel Cardona and Jordi Quer.  \newblock Field of
  moduli and field of definition for curves of genus 2.  \newblock In
  {\em Computational aspects of algebraic curves}, volume~13 of {\em
    Lecture Notes Ser. Comput.}, pages 71--83. World Sci. Publ.,
  Hackensack, NJ, 2005.


\bibitem[CMKT00]{Chao}
Jinhui Chao, Kazuto Matsuo, Hiroto Kawashiro, and Shigeo Tsujii.
\newblock Construction of hyperelliptic curves with {CM} and its application to
  cryptosystems.
\newblock In {\em Advances in cryptology---ASIACRYPT 2000 (Kyoto)}, 
  {\em Lecture Notes in Comput. Sci.}, Vol. 1976, pages 259--273. 
Springer, Berlin, 2000.

\bibitem[CL01]{CL}
Henry Cohn and Kristin Lauter
\newblock {Generating genus 2 curves with complex multiplication}.
\newblock Microsoft Research Internal Technical Report, 2001.


\bibitem[Cox89]{Cox}
David~A. Cox.
\newblock {\em Primes of the form {$x\sp 2 + ny\sp 2$}}.
\newblock John Wiley \& Sons Inc., New York, 1989.

\bibitem[FL06]{FL}
David Freeman and Kristin Lauter.
\newblock Computing endomorphism rings of Jacobians of genus 2 curves.
\newblock Preprint, 2006.

\bibitem[GH00]{GH}
Pierrick Gaudry and Robert Harley.
\newblock Counting points on hyperelliptic curves over finite fields.
\newblock In {\em Algorithmic number theory (Leiden, 2000)},
  {\em Lecture Notes in Comput. Sci.}, Vol. 1838, pages 313--332.
 Springer, Berlin, 2000.

\bibitem[GHKRW]{GHKRW} 
P. Gaudry, T. Houtmann, D. Kohel, C. Ritzenthaler, A. Weng. 
\newblock The 2-adic CM method for genus 2 curves with application to cryptography.
\newblock Advances in Cryptology, ASIACRYPT 2006, Springer-Verlag, LNCS 4284, 114--129, 2006. 



\bibitem[GS05]{GS}
Pierrick Gaudry and {\'E}ric Schost.
\newblock Modular equations for hyperelliptic curves.
\newblock {\em Math. Comp.}, 74(249):429--454 (electronic), 2005.

\bibitem[Gor97]{Goren}
Eyal~Z. Goren.
\newblock On certain reduction problems concerning abelian surfaces.
\newblock {\em Manuscripta Math.}, 94(1):33--43, 1997.

\bibitem[GL04]{GL}
Eyal~Z. Goren and Kristin~E. Lauter.
\newblock {Class invariants for quartic CM fields}. Preprint, 2004.
\newblock to appear in {\em Annales Inst. Fourier}, 2007. 

\bibitem[How95]{Howe}
Everett~W. Howe.
\newblock Principally polarized ordinary abelian varieties over finite fields.
\newblock {\em Trans. Amer. Math. Soc.}, 347(7):2361--2401, 1995.

\bibitem[Igu60]{Igusa1}
Jun-ichi Igusa.
\newblock Arithmetic variety of moduli for genus two.
\newblock {\em Ann. of Math. (2)}, 72:612--649, 1960.

\bibitem[Igu62]{Igusa3}
Jun-ichi Igusa.
\newblock On {S}iegel modular forms of genus two.
\newblock {\em Amer. J. Math.}, 84:175--200, 1962.

\bibitem[Igu67]{Igusa2}
Jun-ichi Igusa.
\newblock Modular forms and projective invariants.
\newblock {\em Amer. J. Math.}, 89:817--855, 1967.

\bibitem[KW89]{KaWa}
Luise-Charlotte Kappe and Bette Warren.
\newblock An elementary test for the {G}alois group of a quartic polynomial.
\newblock {\em Amer. Math. Monthly}, 96(2):133--137, 1989.

\bibitem[Koh96]{Kohel}
David Kohel. 
\newblock Endomorphism rings of elliptic curves over finite fields.
\newblock Ph.D. thesis, University of California, Berkeley, 1996.

\bibitem[Lan73]{Lang}
Serge Lang.
\newblock {\em Elliptic functions}.
\newblock Addison-Wesley Publishing Co., Inc., Reading, Mass.-London-Amsterdam,
  1973.

\bibitem[Mes72]{Messing}
William Messing.
\newblock {\em The crystals associated to {B}arsotti-{T}ate groups: with
  applications to abelian schemes}.
\newblock Springer-Verlag, Berlin, 1972.
\newblock Lecture Notes in Mathematics, Vol. 264.

\bibitem[Mes91]{Mestre}
Jean-Fran{\c{c}}ois Mestre.
\newblock Construction de courbes de genre {$2$} \`a partir de leurs modules.
\newblock In {\em Effective methods in algebraic geometry (Castiglioncello,
  1990)}, volume~94 of {\em Progr. Math.}, pages 313--334. Birkh\"auser Boston,
  Boston, MA, 1991.

\bibitem[Mil98]{Milne}
J.S.\ Milne.
\newblock {\em Abelian Varieties}.
\newblock 1998.
\newblock Available at http://www.\-math.lsa.umich. edu/\~{}jmilne.

\bibitem[OU73]{OortUeno}
Frans Oort and Kenji Ueno.
\newblock Principally polarized abelian varieties of dimension two or three are
  {J}acobian varieties.
\newblock {\em J. Fac. Sci. Univ. Tokyo Sect. IA Math.}, 20:377--381, 1973.


\bibitem[RV00]{RV}
Fernando Rodriguez-Villegas.
\newblock Explicit Models of Genus 2 Curves with Split CM.
\newblock In {\em Algorithmic number theory (Leiden, 2000)}, volume 1838 of
  {\em Lecture Notes in Comput. Sci.}, pages 505--513. Springer, Berlin, 2000.

\bibitem[Shi71]{Shimura71}
Goro Shimura.
\newblock On the zeta-function of an abelian variety with complex
  multiplication.
\newblock {\em Ann. of Math. (2)}, 94:504--533, 1971.

\bibitem[Shi98]{Shimura}
Goro Shimura.
\newblock {\em Abelian varieties with complex multiplication and modular
  functions}, volume~46 of {\em Princeton Mathematical Series}.
\newblock Princeton University Press, Princeton, NJ, 1998.

\bibitem[Spa94]{Spallek}
Anne-Monika Spallek.
\newblock Kurven vom {G}eschlecht 2 und ihre {A}nwendung in
  {P}ublic-{K}ey-{K}ryptosystemen.
\newblock Ph.D. Thesis. Universit\"at Gesamthochschule Essen, 1994.


\bibitem[SW96]{SpeWi}
B.~K. Spearman and K.~S. Williams.
\newblock Relative integral bases for quartic fields over quadratic subfields.
\newblock {\em Acta Math. Hungar.}, 70(3):185--192, 1996.

\bibitem[vW99]{vanWamelen}
Paul van Wamelen.
\newblock Examples of genus two {CM} curves defined over the rationals.
\newblock {\em Math. Comp.}, 68(225):307--320, 1999.

\bibitem[Wen03]{Weng}
Annegret Weng.
\newblock Constructing hyperelliptic curves of genus 2 suitable for
  cryptography.
\newblock {\em Math. Comp.}, 72(241):435--458 (electronic), 2003.

\bibitem[Wen04]{Weng2}
Annegret Weng.
\newblock Extensions and improvements for the CM method for genus two,  
\newblock In: High primes and misdemeanours: lectures in honour of the 60th birthday of Hugh Cowie
Williams, 379--389, Fields Inst. Commun., 41, Amer. Math. Soc., Providence,
RI, 2004.

\end{thebibliography}
\end{document}